\documentclass[11pt]{amsart}
\usepackage{amscd}
\usepackage[arrow,matrix]{xy}
\usepackage{graphicx}
\usepackage{amsmath}
\usepackage{color}

\usepackage{amsmath, latexsym, amssymb}
\input xypic
\numberwithin{equation}{section}
\theoremstyle{plain}
\newtheorem{lemma}{Lemma}[section]
\newtheorem{proposition}[lemma]{Proposition}
\newtheorem{theorem}[lemma]{Theorem}

\theoremstyle{definition}
\newtheorem{definition}[lemma]{Definition}

\begin{document}
\newcommand{\R}{{\mathbb R}}
\newcommand{\C}{{\mathbb C}}
\newcommand{\F}{{\mathbb F}}
\renewcommand{\O}{{\mathbb O}}
\newcommand{\Z}{{\mathbb Z}} 
\newcommand{\N}{{\mathbb N}}
\newcommand{\Q}{{\mathbb Q}}
\renewcommand{\H}{{\mathbb H}}

\newcommand{\Aa}{{\mathcal A}}
\newcommand{\Bb}{{\mathcal B}}
\newcommand{\Cc}{{\mathcal C}}    
\newcommand{\Dd}{{\mathcal D}}
\newcommand{\Ee}{{\mathcal E}}
\newcommand{\Ff}{{\mathcal F}}
\newcommand{\Gg}{{\mathcal G}}    
\newcommand{\Hh}{{\mathcal H}}
\newcommand{\Kk}{{\mathcal K}}
\newcommand{\Ii}{{\mathcal I}}
\newcommand{\Jj}{{\mathcal J}}
\newcommand{\Ll}{{\mathcal L}}    
\newcommand{\Mm}{{\mathcal M}}    
\newcommand{\Nn}{{\mathcal N}}
\newcommand{\Oo}{{\mathcal O}}
\newcommand{\Pp}{{\mathcal P}}
\newcommand{\Qq}{{\mathcal Q}}
\newcommand{\Rr}{{\mathcal R}}
\newcommand{\Ss}{{\mathcal S}}
\newcommand{\Tt}{{\mathcal T}}
\newcommand{\Uu}{{\mathcal U}}
\newcommand{\Vv}{{\mathcal V}}
\newcommand{\Ww}{{\mathcal W}}
\newcommand{\Xx}{{\mathcal X}}
\newcommand{\Yy}{{\mathcal Y}}
\newcommand{\Zz}{{\mathcal Z}}

\newcommand{\zt}{{\tilde z}}
\newcommand{\xt}{{\tilde x}}
\newcommand{\Ht}{\widetilde{H}}
\newcommand{\ut}{{\tilde u}}
\newcommand{\Mt}{{\widetilde M}}
\newcommand{\Llt}{{\widetilde{\mathcal L}}}
\newcommand{\yt}{{\tilde y}}
\newcommand{\vt}{{\tilde v}}
\newcommand{\Ppt}{{\widetilde{\mathcal P}}}
\newcommand{\bp }{{\bar \partial}} 
\newcommand{\ad}{{\rm ad}}
\newcommand{\Om}{{\Omega}}
\newcommand{\om}{{\omega}}
\newcommand{\eps}{{\varepsilon}}
\newcommand{\Di}{{\rm Diff}}

\renewcommand{\a}{{\mathfrak a}}
\renewcommand{\b}{{\mathfrak b}}
\newcommand{\e}{{\mathfrak e}}
\renewcommand{\k}{{\mathfrak k}}
\newcommand{\m}{{\mathfrak m}}
\newcommand{\pg}{{\mathfrak p}}
\newcommand{\g}{{\mathfrak g}}
\newcommand{\gl}{{\mathfrak gl}}
\newcommand{\h}{{\mathfrak h}}
\renewcommand{\l}{{\mathfrak l}}
\newcommand{\sm}{{\mathfrak m}}
\newcommand{\n}{{\mathfrak n}}
\newcommand{\s}{{\mathfrak s}}
\renewcommand{\o}{{\mathfrak o}}
\newcommand{\so}{{\mathfrak so}}
\renewcommand{\u}{{\mathfrak u}}
\newcommand{\su}{{\mathfrak su}}

\newcommand{\ssl}{{\mathfrak sl}}
\newcommand{\ssp}{{\mathfrak sp}}
\renewcommand{\t}{{\mathfrak t }}
\newcommand{\Cinf}{C^{\infty}}
\newcommand{\la}{\langle}
\newcommand{\ra}{\rangle}
\newcommand{\half}{\scriptstyle\frac{1}{2}}
\newcommand{\p}{{\partial}}
\newcommand{\notsub}{\not\subset}
\newcommand{\iI}{{I}}               
\newcommand{\bI}{{\partial I}}      
\newcommand{\LRA}{\Longrightarrow}
\newcommand{\LLA}{\Longleftarrow}
\newcommand{\lra}{\longrightarrow}
\newcommand{\LLR}{\Longleftrightarrow}
\newcommand{\lla}{\longleftarrow}
\newcommand{\INTO}{\hookrightarrow}

\newcommand{\QED}{\hfill$\Box$\medskip}
\newcommand{\UuU}{\Upsilon _{\delta}(H_0) \times \Uu _{\delta} (J_0)}
\newcommand{\bm}{\boldmath}

\title[Realizing homology classes  by  symplectic submanifolds]{\large Realizing homology classes  by  symplectic submanifolds}
\author{ H\^ong V\^an L\^e}
\thanks{partially supported by RVO: 67985840}

\address{Institute of Mathematics of ASCR, Zitna 25, 11567 Praha 1, Czech Republic}
 

\abstract  In this note we prove that a positive multiple of each  even-dimensional integral homology class of a compact symplectic manifold $(M^{2n}, \om)$ can be represented as the difference of the fundamental classes of two  symplectic submanifolds in $(M^{2n}, \om)$.  We also discuss the realizability of integral homology classes by symplectic surfaces in $(M^{2n}, \om)$.
\endabstract

\maketitle
To my Teacher Anatoly Timofeevich Fomenko
\tableofcontents


\section {Introduction.}

In 1954 Ren\'e Thom proved the following celebrated theorem, which relates
the topological structure  with the differentiable structure on compact    smooth manifolds  $M^m$.

\medskip

\begin{theorem}\label{thm:Theorem 1.1.} (\cite[Theorem II.25]{Thom1954}). For each element $\alpha \in H_k (M^m,\Z)$
of  a compact differentiable manifold $M^m$ there exists a positive integer $N(k,m)$
such that the element $N(k,m)\cdot \alpha$ can be realized by a  differentiable submanifold
in $M^m$.
\end{theorem}

Thom's theorem is optimal in the sense that we cannot replace   the positive integer $N(k,m)$ by the number  $ 1$. Namely
Thom showed that for each $k \ge 7$ there is a compact differentiable manifold
$M^m$ and  an element $\alpha \in H_k (M^m, \Z)$ such that $\alpha$ cannot be realized
by the fundamental class of a submanifold in  $M^m$ \cite{Thom1954}.

\medskip

As an immediate consequence of Thom's theorem \ref{thm:Theorem 1.1.} we obtain that  the homology group
$H_*(M^m, \Q)$ is generated by the fundamental classes of differentiable  submanifolds of $M^m$. 
\medskip

In this note we prove the following
weak version of Thom's theorem for compact symplectic manifolds.

 \medskip

\begin{theorem}\label{thm:Theorem 1.2.}  (Main Theorem) Suppose that $(M^{2n},\om)$ is a compact symplectic manifold.  Then  for each element $\alpha\in H_{2k}(M ^{2n}, \Z)$, $1\le k \le n$, there exists
a positive number $N(\alpha)\in \N^+$ such that $N(\alpha) \cdot \alpha = [S_1 ^{2k}]- [S_2 ^{2k}]$,
where $S_1^{2k}$ and $S_2 ^{2k}$ are symplectic submanifolds in $(M^{2n},\om)$. Consequently, $H_{2k}(M^{2n}, \Q)$ is  generated by the fundamental  classes of symplectic submanifolds in $(M^{2n}, \om)$ for all $1\le k \le n$.
\end{theorem}

A particular case of Theorem \ref{thm:Theorem 1.2.} is  the following  celebrated  Theorem due to Donaldson.
For an element $\alpha \in H^*(M, \Z)$ denote by $PD(\alpha)$ the  Poincar\'e dual of  $\alpha$.

\begin{theorem}\label{thm:Theorem 1.3.}(cf. \cite[Theorem 1]{Donaldson1996}) Let $(M^{2n}, \om)$ be a compact integral
symplectic manifold, i.e. the cohomology class $[\om]$ belongs to $H^2 (M^{2n},\Z)
\subset H^2 (M^{2n}, \R)$. Then for each $1 \le k \le n$ there exists a number 
$N_1(k)\in \N^+$ such that $N \cdot PD ([\om ^{n-k}])$ can be realized by 
a symplectic   submanifold $S_1 ^{2k}$ in $(M^{2n}, \om)$, if $N \ge N_1 (k)$.  
\end{theorem}

Actually, in his paper 
Donaldson  stated his result  \cite[Theorem 1] {Donaldson1996}  as  a particular case of Theorem \ref{thm:Theorem 1.3.}  for $k =n- 1$  and using this result  he   derived the existence of symplectic   submanifolds  in other     dimensions  \cite[Corollary 6]{Donaldson1996}. Kaoru Ono
noticed that using a simple  argument we can easily obtain Theorem \ref{thm:Theorem 1.3.} for $1\le k \le n-2$ as a consequence of  Donaldson's results in the cited paper \cite[Proposition 3, Theorem 5]{Donaldson1996}. We shall represent Ono's argument
in section \ref{sec:2}. Theorem \ref{thm:Theorem 1.3.}  also follows from a result by Auroux in  \cite[Theorem 2]{Auroux1997}, see Proposition \ref{rem:2.5} below, whose  proof  generalizes Donaldson's arguments. 

\medskip

We shall call a homology class $\alpha \in H_{2k} (M,\Z)$ {\it symplectic},
if $\alpha$ is the  fundamental class of a symplectic  submanifold in $M^{2n}$.
We shall call a homology class $\alpha \in H_{2k} (M,\Z)$ {\it formal symplectic},
if the pairing $\la \alpha,\om^k\ra $  is positive.


{\bf Conjecture 1.}\label{con:lhv} {\it For each formal symplectic  class $\alpha \in H_{2k} (M^{2n},\Z)$ there exists a positive number $N_1(\alpha)\in \N^+$ such that
$N_1(\alpha)\cdot \alpha$ is symplectic.}



\medskip 

We   also  have the following       

\begin{theorem}\label{thm:Theorem 1.5.} Let $(M^{2n}, \om)$ be a compact symplectic manifold with $2n \ge 6$. 

 a) If $\om|_{\pi _2 (M^{2n})} \not = 0$, then there exists a symplectic sphere in $(M^{2n},
\om)$.

b) Conjecture \ref{con:lhv} is true   for $k =1$.
\end{theorem}

The remainder of this paper is   organized as follows. In section \ref{sec:2} we give a proof of     our Main Theorem.  In section \ref{sec:3} we give  a proof of
Theorem \ref{thm:Theorem 1.5.}. We also  discuss in section \ref{sec:3}   related  results
   concerning    Conjecture \ref{con:lhv}  for $k =1$.

\

This  note  is  the    final version  of   my preprint \cite{Le2005}.

\section {Proof of  the main result}\label{sec:2}

Clearly Theorem  \ref{thm:Theorem 1.2.} is a consequence of  the following 

\medskip

\begin{theorem}\label{thm:2.1}  Let $(M^{2n}, \om)$ be a compact  symplectic manifold
and $\alpha \in H_{2k} (M^{2n}, \Z)$. Then there exist an integral symplectic
form $\bar\om$ on $M^{2n}$ and  positive integral numbers
$N_2(\alpha) $ and $N_3(\alpha)$ such that $[N_2(\alpha) \cdot \alpha  + N_3(\alpha) PD( [\bar\om^{n-k}])$ and $N_3(\alpha)\cdot PD([\bar \om^{n-k}])$ can be realized by  some symplectic submanifolds $S_1 ^{2k}$   and $S_2^{2k}$ of $(M^{2n}, \om)$.
\end{theorem}

The proof of Theorem \ref{thm:2.1} is divided in  two steps.
In the first step we prove  a  particular  case of Theorem \ref{thm:2.1} for  compact integral symplectic  manifolds $(M^{2n}, \om)$, i.e.  we  set $\bar \om:  = \om$. Moreover, for this ``integral" case  the obtained symplectic submanifolds  $S_i ^{2k}$    are  approximately $J$-holomorphic  for some compatible almost complex structure $J$ on $(M^{2n}, \om)$. In the second step, using  a perturbation argument, for an arbitrary symplectic   compact manifold $(M^{2n},\om)$ we shall find a rational symplectic  form $\bar\om$ close
enough to $\om$ such that the   symplectic submanifolds $S_i ^{2k}$ of $(M^{2n}, \bar \om)$, whose existence has been established in the  first step, 
are also  approximately $J$-holomorphic, and hence  they are symplectic submanifolds  of $(M^{2n}, \om)$.

\medskip

The    approximate (resp. asymptotic)  $J$-holomorphy of symplectic  submanfolds  $S^{2k}_i$ are expressed in  terms
of $\eta$-transversality and  approximately  (resp. asymptotically) holomorphic sections, whose zero  sets are $S^{2k}_i$. These notions
have been first introduced by Donaldson for sections of a Hermitian  line bundle over  an almost complex
manifold $(M^{2n},J)$  supplied with a Riemannian metric $g$ \cite{Donaldson1996}.  Then they have been  extended  later by Auroux for
 sections  of a  Hermitian vector bundle over  an almost  complex  manifold $(M^{2n}, J)$ supplied  with a Riemannian metric $g$  \cite{Auroux1997}, \cite{Auroux2000}.

\medskip

\begin{definition}\label{def:2.2} (\cite[Definition 2]{Auroux1997}, cf. \cite[Definition 3.3]{Auroux2000}, \cite[Definition 17]{Donaldson1996}).
Given a constant $\eta > 0$, we say that a section $s$ of a vector bundle $E$ carrying
a metric and a connection $\nabla$ over  a  Riemannian manifold $M^{2n}$ is {\it $\eta$-transverse to $0$}, if at every point $x$ such that
$|s(x)|\le \eta$ the covariant derivative $\nabla s(x): T_xM^{2n}\to E(x)$ is surjective and admits
 a right inverse of norm less than $\eta^{-1}$. Furthermore, we say that  a sequence of sections  of $E$  is {\it uniformly transverse} to $0$  if there exists 
 a fixed constant $\eta$  such that all  sections in the sequence  are $\eta$-transverse to $0$.
\end{definition}

For  a Hermitian  connection $\nabla$  on   a  Hermitian  vector bundle   $E$ over  an almost complex manifold $(M, J)$  let us  denote  by $\bar \p$ the $(0,1)$-part 
(the anti-complex  linear  part) of the covariant derivative $\nabla$.

\begin{definition}\label{def:2.3} (\cite[Definition 1]{Auroux1997}, cf.  \cite[Definition 1]{Auroux2002}, \cite[Definition 3.1]{Auroux2000}, cf. \cite{Donaldson1996}). Let $(M, J, g)$ be an almost-complex
manifold with a Riemannian metric $g$ and let $s_k$ be a  sequence of sections of a  Hermitian vector bundle $E_k$  supplied with a Hermitian connection $\nabla_k$. 
A sequence  of sections 
 $ s_k$ of $E_k$ is   called {\it asymptotically  $J$-holomorphic} with respect to the given connections, if  the following bounds  hold
 $$|s_k| = O(1), \:  |\nabla_k s_k| = O(k ^{1/2}, \: |\bar \p s_k| = o(1),  $$
 $$|\nabla_k \nabla_k s_k| = O(k), \: | \nabla  \bar \p s_k | = O ( k^{1/2}).$$ 
\end{definition}

\medskip

\begin{proof} [Proof of Theorem \ref{thm:2.1}]

\

\underline{Step 1.} Let $(M, \om)$  be a compact integral symplectic  manifold   supplied with a  compatible almost complex structure $J$ and the associated  compatible  metric $g_J$. Let $L$ be a complex line bundle over 
  $(M^{2n},\om)$   such that  $c_1 (L) = [\om]$. Set $L_k : =  L^{\otimes k}$.  First  we  need  the following
 
\begin{proposition}\label{rem:2.5} (\cite[Proposition 1, Theorem 2]{Auroux1997}, cf.  \cite[Corollary 5.2]{Auroux2000}). Let  $E$ be a fixed complex vector bundle
over $(M^{2n}, \om)$. Let $J$  be a compatible  almost  complex structure on $(M^{2n}, \om)$ and $g_J$ the associated Riemannian metric. Then   there  exist  asymptotically $J$-holomorphic sections $s_k$ of
$E\otimes L_k$   which are uniformly transverse to 0 and whose zero sets  $W_k : =s_k ^{-1} (0)$ are smooth symplectic submanifolds in $(M^{2n},\om)$.  Furthermore,  the  submanifold $W_k$
are asymptotically $J$-holomorphic,  i.e. $J(TW_k)$ is within $o(1)$ of  $TW_k$.
\end{proposition}

Another important ingredient   of  of our  proof   is the following  topological fact.

\begin{proposition}\label{prop:2.7}(\cite[Corollary 1.2]{Le2008})   For each homology class $\alpha \in 
H_{2k} (M^{2n}, \Z)$ there exists a number  $N_2 (k,n)\in \N^+$  and a complex vector bundle $E^{n-k}$  of complex dimension $n-k$ over $M^{2n}$  such that
$$ c_{n-k} (E^{n-k}) = N_2 (k, n) \cdot PD(\alpha),$$
$$  c_i (E^{n-k}) = 0, \text{ for all } 1 \le i \le n-k-1.$$
\end{proposition}




\


We shall combine Proposition \ref{rem:2.5}  and  Proposition \ref{prop:2.7}  for the proof of  Theorem \ref{thm:2.1}. For this purpose   we use the following well-known formula for the Chern classes of a tensor
product of complex vector bundles, see e.g. \cite[Appendix A.3]{Hartshorne1977}. Denote by
$c_t (E^r)$ the Chern  polynomial of a  complex vector bundle $E^r$ over $M^{2n}$: $c_t (E^r) =
c_0 (E) + c_1 (E)t + \cdots +c_r (E^r) t^r$. Using the Grothendieck  splitting principle
we write
$$c_t ( E^r) = \Pi _{i =1} ^r (1 + a_i t).$$
Let $F^s$ be another complex vector bundle over $M$ with
$$c_t (F^s) = \Pi _{j =1} ^s (1 + b_j t).$$
Then we have
\begin{equation}
c_t ( E ^r \otimes F^s) = \Pi_{i,j} ( 1 + (a_i + b_j)t).\label{eq:2.9}
\end{equation}
Now let $E^{n-k}$  be  a complex   vector bundle in Proposition \ref{prop:2.7}. It follows from (\ref{eq:2.9}) that   the top Chern class  $c_{n-k }( E^{n-k})$
of $E^{n-k}$  satisfies  the  following   relation  for any $N$
\begin{equation}
c_{n-k }( E^{n-k}\otimes L_N) =  N_2(n-k,n)\cdot PD( \alpha) + N\cdot [\om ^{n-k} ].\label{eq:2.10}
\end{equation}
Now let $N_3(\alpha)$  be  a sufficient  large number   such that there exists a  section $s_{N_3 (\alpha)}$  of  $E^{n-k} \otimes L_{N_3 (\alpha)}$ that  satisfies the condition of Proposition \ref{rem:2.5}.
Taking into account (\ref{eq:2.10})  and Theorem \ref{thm:Theorem 1.3.},   we obtain  Theorem \ref{thm:2.1} for compact integral symplectic manifold $(M^{2n},\om)$, and    thus complete   Step 1.

\medskip

\underline{Step 2.} {\it  Proof of Theorem \ref{thm:2.1}  for a general compact symplectic $(M,\om)$.} We need
the following perturbation result.

\medskip 

\begin{proposition}\label{lem:2.9}  Suppose that $(M^{2n}, \om)$ is a compact symplectic manifold
with a compatible almost complex structure $J$ and  the associated  Riemannian metric $g_J$. Then
 there exist
 an integral symplectic  form $\bar \om$ together with a compatible almost complex structure $\bar J$  and  the associated  compatible Riemannian metric $g_{\bar J}$ over $M^{2n}$ such that the following statement holds.  Let $E^k$ be a Hermitian vector bundle over $(M,\bar \om, \bar J)$ with
a Hermitian connection and $s_k$ be  sections of  $E^k$ which
are uniformly  transverse to $0$ and which are asymptotically $\bar J$-holomorphic. 
Then  for    sufficiently large $k$ the zero section $s^{-1}_k (0)$ is a symplectic submanifold  of $(M^{2n}, \om)$.
\end{proposition}

\begin{proof}  For any symplectic manifold $(M^{2n},\om)$  supplied  with a compatible  almost
complex structure $J$  and  the compatible   metric $g_J$, we denote by $G^\om_{2l}(M^{2n})$ the Grassmannian of $\om$-symplectic
$2l$-planes  in $T_* M^{2n}$ and by $G^J _{2l}(M^{2n})$ the Grassmannian of
$J$-invariant $2l$-planes in $T_*M^{2n}$. Clearly $G^\om _{ 2l} (M^{2n})$ is an open
neighborhood of $G^J_{2l} (M^{2n})$. Denote  by $d_{g_{ J}}(x)$  the     metric on     $\Lambda ^{2l}T_x M^{2n} \supset Gr_{2l} (T_xM^{2n})\supset Gr _{2l} ^\om (T_xM^{2n})$ 
that is induced   by the  Riemannian  metric $g_{ J}(x)$.
The openess  of $G^\om _{ 2l} (M^{2n})$  implies that, if $M^{2n}$ is compact and $\om$ and $J$ are given, there exists a positive number $\eps > 0$
such that  the following  assertion holds  for all $x \in M^{2n}$  and for all $V \in Gr_{2l} (T_x M^{2n})$
\begin{equation}
d_{g_{ J}}( V^{2l}(x), G^J_{2l} (T_xM^{2n})) < \eps  \LRA  V^{2l} \in  Gr_{2l}^\om(T_xM^{2n}).\label{eq:jiom}
\end{equation} 

\begin{lemma}\label{lem:per} 
For all $ l \in [1, n-1]$  there  exists  a small  number  $\eps$ such that   the following assertion 
\begin{equation}
d_{g_{J_1}} (V, G^{J_1}_{2l} (T_x M^{2n})) < (\eps /2)  \LRA V \in   Gr^{\om}_{2l}(T_xM)\label{eq:both}
\end{equation}
holds for   all $x \in  M^{2n}$  for all $V \in Gr_{2l} (T_x M^{2n})$ and   for  all    almost  complex  structure $J_1$   on $M^{2n}$ that satisfies  the following condition. $J_1$  is  compatible  with  some  symplectic  structure $\om_1$  such that 
\begin{equation}
|\om -\om_1| _{g_J}  +|J-J_1|_{g_J}  < {\eps\over 16}.\label{eq:eps}
\end{equation}
\end{lemma}
\begin{proof}   There  exists   a number $0< \eps  < 1/2$ such   if (\ref{eq:eps}) holds  then
\begin{equation}
d_{g_{J_1}} ( V,  Gr^{J_1} _{2l} (T_x M^{2n})) < {\eps \over 2}  \LRA   d _{ g_{J_1}} (V, Gr ^{J}_{2l} (T_x M^{2n}) )<  {\eps \over 2} + {\eps \over 8}\label{eq:eps1}
\end{equation}
for all $x \in M$ and for all $V \in Gr _{2l} (T_x M^{2n})$. 
Taking  into account (\ref{eq:eps}), and shrinking   $\eps$ if necessary,   we obtain from (\ref{eq:eps1})
\begin{equation}
d_{g_{J}} (V, Gr ^{J}_{2l} (T_x M^{2n})) < {3\over 2} d_{ g_{J_1}} (V, Gr ^{J} (T_x M^{2n})) < {3 \eps \over 4} + { 3\eps  \over 16}.\label{eq:eps2}
\end{equation}
Taking into account (\ref{eq:jiom}) we obtain    (\ref{eq:both})  immediately from (\ref{eq:eps2}). This completes  the proof of  Lemma \ref{lem:per}.
\end{proof}
Now let us complete  the    proof  of Proposition  \ref{lem:2.9}. Let $\eps$ be  a  small positive number  in Lemma \ref{lem:per}.  Let  $\om_1$    a rational  symplectic form  on $(M^{2n}, \om)$  such that    there exists  a compatible almost complex structure $J_1$   for  which   the  condition (\ref{eq:eps}) holds. Suppose that   $\bar \om$ is a multiple of $\om_1$ such that
$[\bar \om]\in H^2 (M, \Z)$. Then the compatible almost
complex structure $J_1$ is also compatible to $\bar \om$. Now we set $\bar J : = J_1$.
Next we note that the metrics $g_1$ and $\bar g$ associated respectively  to $(\om_1, J_1)$ and
$(\bar \om, \bar J)$ induce the same  metric $d_{g_1}$ on  each $G_{2l} (T_x M)$. Furthermore,
any $ \om_1$-symplectic plane (resp. $J_1$-invariant plane)  is also $\bar \om$-symplectic (resp. $\bar J$-invariant plane) and  the  converse  also holds.

Now let $s_k$ be  sections of  a complex vector bundle  $E^k$  supplied   with a  Hermitian connection    as in Proposition \ref{lem:2.9}   with respect to $\bar J$.  By Proposition \ref{rem:2.5},  for  sufficiently large $k$,  the   zero  sections $W_k:=s_k ^{-1} (0)$  are asymptotically $\bar J$-holomorphic.  Lemma \ref{lem:per} implies  that
$W_k$  are   symplectic  submanifolds  of $(M, \om)$ and  of $(M, \bar \om)$  for sufficiently large  $k$. This  completes  the proof of Proposition \ref{lem:2.9}.
\end{proof}

{\it Completion of the proof of Theorem \ref{thm:2.1}.} Now given  a compact symplectic
manifold $(M,\om)$ with a compatible almost complex structure $J$  and $\alpha\in
H_2 (M,\Z)$ we shall choose $\bar \om$ and $\bar J$ as in Proposition \ref{lem:2.9}. Denote by $L_{\bar \om}$  the  vector line bundle  with $c_1(L) = \bar \om$. Let $E^{n-k}$ be a  complex vector bundle  satisfies the condition of Proposition \ref{prop:2.7}.  We set $E_p^{n-k}: = E^{n-k}\otimes 
L^{\otimes p}_{\bar \om}$.  Let $s_p$ be  sequence of  asymptotically $\bar J$-holomorphic sections of $E_p^{n-k}$  whose existence  is ensured by  Proposition \ref{rem:2.5}.  It follows from
Proposition \ref{lem:2.9} that the zero sections  $s^{-1}_p(0)$ are symplectic submanifolds of $(M^{2n},\om)$, for sufficiently large $p$. Taking into account   (\ref{eq:2.10})  this completes  the proof  of Theorem \ref{thm:2.1}.
\end{proof}

\medskip

The remainder of this section is devoted to  a simple  proof of Theorem \ref{thm:Theorem 1.3.}   due to Kaoru Ono.

\begin{proof}[An elementary proof of Theorem \ref{thm:Theorem 1.3.} for $1\le k \le n-2$] \cite{Ono2004}.
We note that Theorem \ref{thm:Theorem 1.3.} is a consequence
of the following 

\begin{proposition}\label{prop:ono} There are positive integral numbers $n_1, \cdots, n_k$
and  for each $i = \overline{1,k}$ a section $s_i$ of the line bundle  $L^{N_i}_\om$ such that the section $\hat s: =s_1 \oplus  \cdots \oplus s_k$
intersects to the zero section of $\hat L : =L^{N_1}_\om \oplus \cdots \oplus L^{N_k}_\om$ transversally and $\hat s^{-1} (0)$ is a symplectic submanifold in $(M^{2n}, \om)$,
if $N_i \ge n_i$ for $ 1\le i \le k$.
\end{proposition}
\begin{proof}
We prove Proposition \ref{prop:ono} inductively on $k$.  For
$k =1$ the statement  is exactly Theorem 1 in \cite{Donaldson1996}.
Assume that the above statement is valid for $k = K$. We shall prove its validity for 
$k = K +1$. We denote by $S_K$ the common zero locus of sections $s_1, \cdots, s_K$
which is by the induction assumption a symplectic submanifold. We note that
the restriction of the line bundle $L_\om$ to $S_K$ has the curvature $\om_{| S_K}$ which
is the symplectic form on $S_K$.
 According to  Donaldson  \cite[Theorem 5 and Proposition 3]{Donaldson1996} there exists a number
 $N_0$ such that if  $n_{K+1} > N_0$ there is a section 
 $$ \tilde s_{K+1} : S_K\to L^{n_{K+1}}_\om$$ such that  
$\tilde s_{K+1}$ intersects with the zero section of $L^{n_{K+1}}_\om$ transversally.
Moreover  the zero section  $N:=\tilde s_{K+1}^{-1}(0) $  is a symplectic submanifold
in $S_K$.

Since the fiber $L^{n_{K+1}}_\om$ is contractible, the section $\tilde s_{K +1}$ can be extended to a section $s_{K+1} : M ^{2n} \to L^{n_{K+1}}_\om$. The proof of the above statement is complete, if we can show
that $s_{K+1}$ can be chosen to be  transversal to   the zero section of
a neighborhood $U_\eps(S_K)\subset M^{2n}$.
It is easy to see that the extension of the section from
the submanifold to the ambient space is just done by
pull back the section by the projection of the normal bundle of the symplectic 
submanifold. More explicitly,  we choose $U_\eps (S_K)$  to be a (geodesic) neighborhood of $S_K$ in $M^{2n}$,
such that there is a diffeomorphism $f$ from $U_\eps (S_K)$ to a neighborhood
$V_\eps (S_K)$ of the zero
section of the  normal bundle $V(S_K)$ of $S_K$ in $M$. By using this diffeomorphism
$f$ we can work now on $V_\eps (S_K)$. It is easy to see  that 
$$V(S_K) = (L^{n_1}_\om\oplus \cdots \oplus L^{n_K}_\om)_{|S_K}.$$
Then we let the extension of $\tilde s_{K+1} : S_K \to  L^{n_{K+1}}$ to a section
$ \bar s_{K+1} : V_\eps (S_K) \to (f^{-1}) ^* L^{n_{K+1}}$ be defined by
$$ \bar  s_{K+1} (x,l_1, \cdots, l_{K} )  = ( x,l_1, \dots , l_K, \tilde s_{K+1} (x) ).$$
Finally we extend $\bar s_{K+1}$ to a section $s_{K+1}$ over $M^{2n}$. This section
$s_{K+1}$  satisfies the condition of  Proposition \ref{prop:ono}.
\end{proof}
This completes the proof of Theorem \ref{thm:Theorem 1.3.}
\end{proof}
We  observe that Proposition \ref{prop:ono} is  also a  consequence of Proposition \ref{rem:2.5}.

\section {Realizing  homology classes  by symplectic surfaces}\label{sec:3}

In this  section we  give  a  proof  of  Theorem \ref{thm:Theorem 1.5.}  and discuss  results  related to  Conjecture \ref{con:lhv}.

\begin{proof}[Proof of Theorem \ref{thm:Theorem 1.5.}] To prove  the first assertion of  Theorem \ref{thm:Theorem 1.5.} we need the following

\begin{lemma}\label{lem:3.1} Suppose that the condition  of Theorem \ref{thm:Theorem 1.5.}.a is satisfied. Let us  denote by $\om_0$ the standard symplectic form on $S^2$.Then
there is an embedding $f:S^2\to M^{2n}$ such that $f^*([\om]) = k[\om_0]$ for some $k> 0$. 
\end{lemma}

\begin{proof} We provide    $M^{2n}$ with a  Riemannian metric $g$. We  claim that   the required  embedding  can be  obtained by  a perturbation of
a minimal immersion of  $S^2 \to (M^{2n}, g)$  whose  existence  follows from the following theorem due to Sacks and Uhlenbeck

\begin{proposition}\label{prop:su} (\cite[Theorem 5.9]{SU1981}) There exists  a set of free homotopy classes $\Lambda_i \in \pi_0 C^0(S^2, M^{2n})$ such that elements $\lambda \in \Lambda_i$  generate $\pi_2 (M^{2n})$ acted  on by $\pi_1 (M^{2n})$ and  each $\Lambda_i$ contains  a conformal  branch immersion  of a sphere   having least area  among maps of $S^2$ into $N$  which lies  in $\Lambda_i$.
\end{proposition} 
 Since $\om |_{\pi_2  (M^{2n}) }\not  = 0$ it follows  that  at least  one  of  the minimal immersions $f$  obtained by  Proposition \ref{prop:su}
 satisfies  the   condition $f ^*[(\om)]  = k [\om]$ for some  $ k > 0$ (we can change the orientation of the map $f$  if $k < 0$).
 Since $f$ is a  branch minimal immersion having least area,  it has only isolated  singular  points. Now, using  isotopy  and  dimension  condition $2n \ge 6$ we     perturb $f$ slightly  to  obtain  the required embedding.
\end{proof}

\medskip

{\it Continuation of the proof of Theorem \ref{thm:Theorem 1.5.}.a} First we  find  a map covering  $F : T_* S^2 \to  T_* M^{2m}$ of $f$ which is fiber-wise
symplectic.  In other words  we     find a section  of the bundle $Iso_{sym} (T_*S^2, T_*M^{2m})$ over $S^2$
whose fiber  is $Symp(2m) / Symp (2m-2)$.  The fiber of this bundle
  which is $2m-1$ connected  \cite{Gromov1986}. So there is no obstruction for such  a section.
Now Theorem \ref{thm:Theorem 1.5.} follows immediately from the Gromov h-principle, stated  below,  and from the observation, that we can make a $C^1$-perturbation of a symplectic immersion  of
$S^2$ to get a symplectic embedding, since the dimension of $M^{2n}$ is at least $6$.

\medskip

{\bf H-principle for symplectic immersions} \cite[3.4.2.A]{Gromov1986}  {\it Let $\om_0$ be 
an arbitrary  (possibly singular) closed smooth 2-form on a smooth manifold $V$
and let $F_0 : (T_*V, \om_0) \to (T_*M, \om)$ be a fiberwise injective isometric homomorphism. Let us denote by $i$ the embedding $V\to T^*V$ as the zero section, and by $\pi$ the projection  $T^* M \to M$.
Assume that $(\pi \circ F \circ i) ^* [\om] = [\om_0]$.  
If $\dim V  < \dim M$, then the map
$\pi \circ F \circ i: V \to M ^{2n}$ admits a fine $C^0$- approximation by isometric  smooth immersions
$f: (V,\om_0) \to (M^{2n}, \om)$ whose differentials $Df: T_*V \to T_*M$ are homotopic to $F_0$ in the space of fiberwise injective
isometric homomorphisms.}

This completes  the  proof  of Theorem \ref{thm:Theorem 1.5.}.a.

To  prove the second  assertion  of Theorem \ref{thm:Theorem 1.5.}  we use   Theorem \ref{thm:Theorem 1.1.} due to Thom instead  of Lemma \ref{lem:3.1}.
The rest  of the argument  of the   proof repeats  those  in the  proof  of the first  assertion, so we omit it here.
This completes  the  proof   of Theorem \ref{thm:Theorem 1.5.}.
\end{proof}

In \cite{Li} Li  gives   a proof  of the following Theorem, which  strengthens    the second assertion of  Theorem \ref{thm:Theorem 1.5.}

\begin{proposition}\label{prop:li}(\cite[Theorem 1]{Li})  Suppose $(M^{2n}, \om)$ is a symplectic manifold of dimension $2n$. Let
$A$ be any $\om$-positive class in $H_2(M; \Z)$, i.e. $A$ is formal symplectic. Then
\begin{enumerate}
\item $A$ is represented by a connected embedded $\om$-symplectic surface if $2n\ge 6$.
\item A is represented by a connected immersed $\om$-symplectic surface if $2n= 4$.
\end{enumerate}
\end{proposition}

Recently  in \cite{Hamilton} Hamilton  disproved  Conjecture \ref{con:lhv}  for    dimension $2n =4$  by presenting  an obstruction for  presenting
a homology class  of a  symplectic 4-manifold  by an embedded  symplectic surface.

\section*{Acknowledgement} I am  thankful to  Simon
Donaldson, J\"urgen Jost   and Pavel Pudlak for their interests and
supports. I am  indebted to Kaoru Ono for his 
help  and criticism. I am grateful to Dietmar Salamon for pointing out an error in an earlier version of this note and for  encouraging  me to submit the final version  of this note.      I thank  Tian-Jun Li   for sending   me  his paper \cite{Li}. I   acknowledge  MFI   at ETH in Z\"urich  for  their hospitality and financial support in 2005,    where   a part of this paper has been written.
 It is my pleasure to dedicate this note to my teacher Anatoly Timofeevich Fomenko, for he introduced me into  algebraic topology and symplectic geometry and for his support and encouragement. 

\end{document}